\documentclass[12pt]{article}
\usepackage{amssymb}
\usepackage[margin=1in]{geometry}
\usepackage{mathrsfs}
\usepackage{stmaryrd}
\usepackage{amsfonts,amsmath,amssymb,amscd,texdraw, colordvi}
\usepackage{shadow}
\usepackage{graphicx}
\usepackage{color}
\usepackage{longtable}
\usepackage{pstricks,multido}
\usepackage{hyperref}
\allowdisplaybreaks

\newtheorem{thm}{Theorem}[section]

\newtheorem{lem}[thm]{Lemma}

\newtheorem{conj}{Conjecture}[section]

\def\qed{\hfill \rule{4pt}{7pt}}

\def\pf{\noindent {\it{Proof.} \hskip 2pt}}
\numberwithin{equation}{section}

\begin{document}
\begin{center}
{\large\bf  Proof of a conjecture of  Morales--Pak--Panova on

reverse plane partitions}
\end{center}

\begin{center}

{\small Peter L. Guo$^1$,   C.D. Zhao$^2$ and Michael X.X. Zhong$^3$}

\vskip 4mm
$^{1,2}$Center for Combinatorics, LPMC-TJKLC\\
Nankai University,
Tianjin 300071,
P.R. China \\[3mm]

$^3$College of Science\\
Tianjin University of Technology\\
Tianjin 300384, P.R. China

\vskip 4mm

$^1$lguo@nankai.edu.cn\\
$^2$2120150006@mail.nankai.edu.cn, $^3$zhong.m@tjut.edu.cn

\end{center}

\begin{abstract}
Using  equivariant  cohomology theory,
Naruse obtained  a hook length formula for the number of standard Young tableaux of  skew shape $\lambda/\mu$. Morales,
Pak and Panova  found  two $q$-analogues of Naruse's formula
respectively by counting  semistandard Young tableaux of shape $\lambda/\mu$ and reverse plane partitions of shape $\lambda/\mu$.
When $\lambda$ and $\mu$ are  both staircase shape partitions,
 Morales, Pak and Panova  conjectured that the generating function of reverse plane partitions of shape  $\lambda/ \mu$ can   be expressed as a determinant whose
 entries are   related to $q$-analogues of the Euler numbers.
 The objective of this paper is to prove this conjecture.
\end{abstract}



\section{Introduction}

In the context of equivariant Schubert calculus,
Naruse \cite{Nar} presented a hook length formula for the number of standard Young tableaux of  skew shape $\lambda/\mu$. Recently, Morales,
Pak and Panova \cite{MPP} provided  two $q$-analogues of Naruse's formula
by considering the generating function of  semistandard Young tableaux  of shape $\lambda/\mu$, as well as the generating function of reverse plane partitions of shape $\lambda/\mu$. Denote by $\delta_n=(n-1,n-2,\ldots,1)$   the staircase shape partition.
In the case when $\lambda=\delta_{n+2k}$ and $\mu=\delta_n$,
Morales, Pak and Panova \cite{MPP} conjectured that  the generating function of reverse plane partitions of shape $\lambda/\mu$  can also be expressed as a determinant with entries  determined by $q$-analogues of the Euler numbers. They confirmed   the conjecture  for $k=1$. In this paper, we prove  that   this conjecture is true  for any positive integer $k$. After the completion of this paper, we noticed that Hwang, Kim,  Yoo and Yun \cite{Kim} independently proved
this conjecture via a different approach.

\begin{figure}[h t]
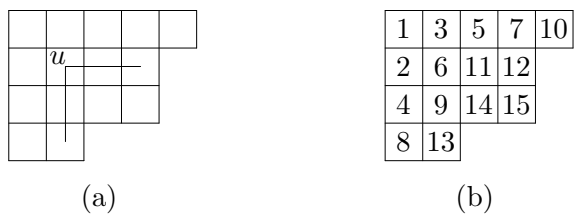

\centertexdraw{ \drawdim mm \linewd 0.1 \setgray 0

\move(0 5)\lvec(10 5) \move(0 10)\lvec(20 10)
\move(0 15)\lvec(20 15) \move(0 20)\lvec(25 20)
\move(0 25)\lvec(25 25)

\move(0 5)\lvec(0 25)\move(5 5)\lvec(5 25)
\move(10 5)\lvec(10 25)\move(15 10)\lvec(15 25)
\move(20 10)\lvec(20 25)\move(25 20)\lvec(25 25)

\move(7.5 17.5)\lvec(17.5 17.5)\move(7.5 17.5)\lvec(7.5 7.5)
\textref h:C v:C \htext(6.5 18.7){{\small $u$}}

\move(50 5)\lvec(60 5) \move(50 10)\lvec(70 10)
\move(50 15)\lvec(70 15) \move(50 20)\lvec(75 20)
\move(50 25)\lvec(75 25)

\move(50 5)\lvec(50 25)\move(55 5)\lvec(55 25)
\move(60 5)\lvec(60 25)\move(65 10)\lvec(65 25)
\move(70 10)\lvec(70 25)\move(75 20)\lvec(75 25)

\textref h:C v:C \htext(52.5 22.5){{\small 1}}
\textref h:C v:C \htext(57.5 22.5){{\small 3}}
\textref h:C v:C \htext(62.5 22.5){{\small 5}}
\textref h:C v:C \htext(67.5 22.5){{\small 7}}
\textref h:C v:C \htext(72.5 22.5){{\small 10}}
\textref h:C v:C \htext(52.5 17.5){{\small 2}}
\textref h:C v:C \htext(57.5 17.5){{\small 6}}
\textref h:C v:C \htext(62.5 17.5){{\small 11}}
\textref h:C v:C \htext(67.5 17.5){{\small 12}}
\textref h:C v:C \htext(52.5 12.5){{\small 4}}
\textref h:C v:C \htext(57.5 12.5){{\small 9}}
\textref h:C v:C \htext(62.5 12.5){{\small 14}}
\textref h:C v:C \htext(67.5 12.5){{\small 15}}
\textref h:C v:C \htext(52.5 7.5){{\small 8}}
\textref h:C v:C \htext(57.5 7.5){{\small 13}}

\textref h:C v:C \htext(12 0){{\small (a)}}
\textref h:C v:C \htext(62 0){{\small (b)}}
} \caption{A Young diagram  and a standard Young tableau}\label{Y-1}
\end{figure}

The classical hook length formula due to Frame Robinson and  Thrall \cite{F-R-T} gives a product formula for the number of
standard Young tableaux whose shape is a Young diagram.  Let $\lambda=(\lambda_1,\lambda_2,\ldots,\lambda_\ell)$ be a partition of a nonnegative  integer $n$, that is, $\lambda=(\lambda_1,\lambda_2,\ldots,\lambda_\ell)$ is a sequence of nonnegative  integers such that $\lambda_1\geq \lambda_2\geq \cdots \geq \lambda_\ell\geq 0 $ and $\lambda_1+\lambda_2+\cdots+\lambda_\ell=n$.
The Young diagram
of  $\lambda$ is a left-justified array of squares  with $\lambda_i$ squares in  row $i$. We   use $(i,j)$ to represent  the square   in row $i$ and column $j$. The hook length $h_u$ of a square $u$ is  the number of squares directly to the right   or directly below $u$, counting $u$ itself once.
For example,  Figure \ref{Y-1}(a)  illustrates the Young diagram of  $(5,4,4,2)$ where the hook  length of $u=(2,2)$ is $5$.

A standard Young tableau of shape $\lambda$ is an assignment of positive integers
$1,2,\ldots,n$ to the squares of $\lambda$ such that the numbers are increasing in
each row and in each column, see Figure \ref{Y-1}(b) for an example.
Let $f^\lambda$ be the number of standard Young tableaux of shape $\lambda$.
The hook length formula \cite{F-R-T} states that
\begin{equation}\label{Hfor}
f^\lambda=\frac{n!}{\prod_{u\in \lambda}h_u}.
\end{equation}

\begin{figure}[h t]
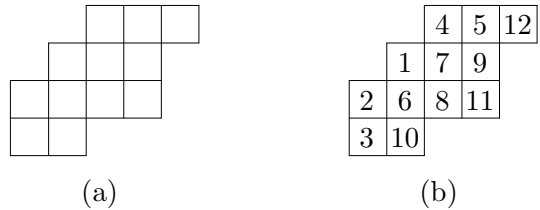

\centertexdraw{ \drawdim mm \linewd 0.1 \setgray 0

\move(45 5)\lvec(55 5) \move(45 10)\lvec(65 10)
\move(45 15)\lvec(65 15) \move(50 20)\lvec(70 20)
\move(55 25)\lvec(70 25)

\move(45 5)\lvec(45 15)\move(50 5)\lvec(50 20)
\move(55 5)\lvec(55 25)\move(60 10)\lvec(60 25)
\move(65 10)\lvec(65 25)\move(70 20)\lvec(70 25)

\move(90 5)\lvec(100 5) \move(90 10)\lvec(110 10)
\move(90 15)\lvec(110 15) \move(95 20)\lvec(115 20)
\move(100 25)\lvec(115 25)

\move(90 5)\lvec(90 15)\move(95 5)\lvec(95 20)
\move(100 5)\lvec(100 25)\move(105 10)\lvec(105 25)
\move(110 10)\lvec(110 25)\move(115 20)\lvec(115 25)

\textref h:C v:C \htext(57 0){{\small (a)}}
\textref h:C v:C \htext(102 0){{\small (b)}}

\textref h:C v:C \htext(92.5 7.4){{\small 3}}
\textref h:C v:C \htext(97.5 7.4){{\small 10}}
\textref h:C v:C \htext(92.5 12.4){{\small 2}}
\textref h:C v:C \htext(97.5 12.4){{\small 6}}
\textref h:C v:C \htext(102.5 12.4){{\small 8}}
\textref h:C v:C \htext(107.5 12.4){{\small 11}}
\textref h:C v:C \htext(97.5 17.4){{\small 1}}
\textref h:C v:C \htext(102.5 17.4){{\small 7}}
\textref h:C v:C \htext(107.5 17.4){{\small 9}}
\textref h:C v:C \htext(102.5 22.4){{\small 4}}
\textref h:C v:C \htext(107.5 22.4){{\small 5}}
\textref h:C v:C \htext(112.5 22.4){{\small 12}}

} \caption{A skew diagram  and a standard Young tableau of shew shape}\label{RPP}
\end{figure}

A standard Young tableau of  skew shape is defined similarly. Recall that a skew diagram  $\lambda/\mu$
is obtained from $\lambda$ by removing the squares of $\mu$, where $\mu$ is a partition whose Young diagram is contained in $\lambda$. For example,
Figure \ref{RPP}(a) is   the skew diagram $\lambda/\mu$ with
 $(5,4,4,2)$ and $\mu=(2,1)$. We use  $|\lambda/\mu|$ to represent  the number of
 squares in $\lambda/\mu$. A standard Young tableau of shape $\lambda/\mu$ is  an
assignment of positive integers
$1,2,\ldots, |\lambda/\mu|$ to the squares of $\lambda/\mu$ such that the numbers are increasing in
each row and  each column, see   Figure \ref{RPP}(b) for an example.

The formula of Naruse \cite{Nar} for the number of standard Young tableaux of shape $\lambda/\mu$ can be expressed as a sum of products  of  hook lengths over excited diagrams of $\lambda/\mu$. The structure of   an excited diagram of $\lambda/\mu$ was introduced independently by Ikeda and Naruse \cite{I-N}, Kreiman \cite{Kre-1,Kre-2}, and Knutson,  Miller and  Yong \cite{K-M-Y}.
 To define  an excited diagram, let us recall the operation  of an excited move on squares of  $\lambda$.
 Let $S$ be a subset of the squares  of $\lambda$. A square $(i,j)$ in $S$
 is called active  if the three squares $(i,j+1), (i+1,j), (i+1,j+1)$ belong to
$\lambda$ but not belong to $S$. An excited move on an active square $(i,j)$ in $S$ is a replacement of
the active square  $(i,j)$ by the square $(i+1,j+1)$. An excited diagram of $\lambda/\mu$ is a
subdiagram  of $\lambda$ obtained from the Young diagram of $\mu$ by a sequence of
excited moves on active squares. Let $\mathrm{E}(\lambda/\mu)$  be the set of
excited diagrams of $\lambda/\mu$. Let $f^{\lambda/\mu}$ be the number of standard Young tableaux of shape $\lambda/\mu$.  Naruse's formula states that
\begin{equation}\label{Nfor}
f^{\lambda/\mu}=|\lambda/\mu|!\sum_{D\in \mathrm{E}(\lambda/\mu)}\prod_{u\in \lambda\setminus D} \frac{1}{h_u}.
\end{equation}
Notice  that when $\mu$ is the empty shape,    \eqref{Nfor} reduces to the  hook length formula \eqref{Hfor} for an ordinary shape.

Morales,
Pak and Panova \cite{MPP} proved two   $q$-analogues of Naruse's formula
respectively by counting semistandard Young tableaux of shape $\lambda/\mu$ and reverse plane partitions of shape $\lambda/\mu$. Recall that a semistandard Young tableau of shape $\lambda/\mu$ is an assignment of positive   integers to the squares of $\lambda/\mu$ such that the numbers are weakly increasing in each row and strictly increasing in each column.
Let $\mathrm{SSYT}(\lambda/\mu)$ be the set of semistandard Young tableaux of
shape $\lambda/\mu$.
Based on the  properties of factorial Schur functions due to Ikeda and Naruse \cite{I-N}, Knutson and Tao \cite{K-T}
and Lakshmibai,  Raghavan and  Sankaran \cite{L-R-S}, Morales,
Pak and Panova \cite{MPP} deduced  the following $q$-analogue of formula \eqref{Nfor}:
\begin{equation}\label{formula-skew}
\sum_{T\in \mathrm{SSYT}(\lambda/\mu)} q^{|T|}=\sum_{D\in \mathrm{E}(\lambda/\mu)} \prod_{(i,j)\in \lambda\setminus D}\frac{q^{\lambda_j'-i}}{1-q^{h_{(i,j)}}},
\end{equation}
where $|T|$ denotes the sum of entries  in $T$, and $\lambda_j'$ is the number of squares of $\lambda$ in column $j$.
As explained in \cite{MPP}, combining the tool of $P$-partitions of Stanley, Naruse's formula \eqref{Nfor} is  a
direct consequence of \eqref{formula-skew}.

The second  $q$-analogue of Naruse's formula is  obtained by   counting reverse plane partitions of shape $\lambda/\mu$.
A   reverse plane partition $\pi$ of shape $\lambda/\mu$ is an assignment of nonnegative integers to the squares of $\lambda/\mu$ such that the numbers are weakly increasing in each row and each column.
Morales, Pak and Panova \cite{MPP} showed that the generating function
of reverse plane partitions of shape $\lambda/\mu$ can be expressed
as a sum over pleasant diagrams of $\lambda/\mu$.
A pleasant diagram of $\lambda/\mu$ is a subset $S$ of squares of $\lambda$
such that $S\subseteq \lambda\setminus D$ for some active diagram $D$ of $\lambda/\mu$.  Let $\mathrm{RPP}(\lambda/\mu)$ be  the set of reverse
plane partitions of shape $\lambda/\mu$, and let $\mathrm{P}(\lambda/\mu)$
 be the set of pleasant diagrams of $\lambda/\mu$.
Utilizing the properties of the Hillman-Grassl algorithm, Morales, Pak and Panova \cite{MPP} deduced that
\begin{equation}\label{GRP}
\sum_{\pi\in \mathrm{RPP}(\lambda/\mu)}{q^{|\pi|}}=\sum_{S\in \mathrm{P}(\lambda/\mu)} \prod_{u\in S} \frac{q^{h_u}}{1-q^{h_u}},
\end{equation}
where  $|\pi|$ denotes the sum of entries of $\pi$.
Similarly, together with Stanley's  $P$-partition technique, \eqref{GRP}
 yields Naruse's formula  \eqref{Nfor}.

Morales, Pak and Panova \cite{MPP} conjectured that when $\lambda$ and $\delta$ are staircase shapes, the generating function in \eqref{GRP} can also be expressed as the determinant with entries determined by $q$-analogues of the Euler numbers.
The   Euler numbers $E_n$ are positive integers defined  by
 \[\sum_{n\geq 0} E_n\frac{x^n}{n!}=\sec x +\tan x,\]
 see, for example, the survey \cite{Sta} of Stanley.
It is  known that $E_n$ counts the number of alternating permutations or reverse alternating permutations on
the set  $[n]=\{1,2,\ldots,n\}$. A permutation $\sigma=\sigma_1\sigma_2\cdots \sigma_n$ on $[n]$ is called an alternating  permutation if
\[\sigma_1> \sigma_2< \sigma_3>\sigma_4< \cdots ,\]
and a reverse  alternating  permutation   if
\[\sigma_1< \sigma_2> \sigma_3<\sigma_4> \cdots .\]
Clearly, alternating permutations on $[n]$ have the
same number as reverse alternating permutations on $[n]$.

Morales, Pak and Panova \cite{MPP} defined a  $q$-analogue $E^{*}_{2n+1}(q)$   of  $E_{2n+1}$ based on reverse alternating permutations.
For a permutation
 $\sigma$ on $[n]$, let $\mathrm{maj}(\sigma)$ be  the major index of $\sigma$, namely,
 \[\mathrm{maj}(\sigma)=\sum_{i\in \mathrm{Des(\sigma)}}i,\]
 where $\mathrm{Des}(\sigma)=\{i\,|\,1\leq i\leq n-1, \sigma_i>\sigma_{i+1}\}$ is
 the  set of descents of $\sigma$.
Set
\[E^{*}_{2n+1}(q)=\sum_{\sigma\in \mathrm{Ralt}_{2n+1}} q^{\mathrm{maj}(\sigma^{-1}\kappa )},\]
where
$\mathrm{Ralt}_{2n+1}$ is the set of reverse alternating permutations on $[2n+1]$, and $\kappa$ is the reverse alternating  permutation $132\cdots(2n+1)(2n)$.
Here,  for a reverse alternating permutation  $\sigma\in \mathrm{Ralt}_{2n+1}$,
  $\sigma^{-1}\kappa$ is   obtained from $\sigma^{-1}$ by interchanging the values $2i$ and $2i+1$ for $1\leq i\leq n$.

Let $\delta_n$ represent the staircase shape partition $(n-1,n-2,\ldots,1)$.
Morales, Pak and Panova \cite{MPP} posed the following conjecture on the
generating function of revers plane partitions of shew shape $\delta_{n+2k}/\delta_n$.

\begin{conj}[\mdseries{Morales, Pak and Panova \cite{MPP}}]\label{conjecture}
For any positive integers $n$ and $k$,
\begin{equation}\label{mainresult}
q^{N}\sum_{\pi\in \mathrm{RPP}(\delta_{n+2k}/\delta_n)}{q^{|\pi|}}
=\det\left[\widetilde{E}^{*}_{2(n+i+j)-3}(q)\right]_{i,j=1}^{k},
\end{equation}
where
\[N=\frac{k(k-1)(6n+8k-1)}{6},\]
and for any positive integer $m$,
\[\widetilde{E}^{*}_{2m+1}(q)=\frac{E^{*}_{2m+1}(q)}{(1-q)\cdots(1-q^{2m+1})}.\]
\end{conj}

Using $P$-partitions,   Morales, Pak and Panova \cite{MPP}
proved   Conjecture \ref{conjecture}   for the case $k=1$, that is,
\begin{equation}\label{k=1}
\sum_{\pi\in \mathrm{RPP}(\delta_{n+2}/\delta_n)}{q^{|\pi|}}=\widetilde{E}^{*}_{2n+1}(q).
\end{equation}

In this paper, we aim to prove Conjecture \ref{conjecture}.
We introduce the structure of an staircase alternating array
 which is an array of alternating words subject to certain conditions.
 We show that  the
 determinant  in Conjecture \ref{conjecture} is the (signed) weight generating function of staircase alternating arrays.
 We construct  an involution $\Phi$ on staircase alternating arrays.
As a consequence,  the determinant in Conjecture \ref{conjecture} is a weighted
counting of   the fixed points of $\Phi$. On the other hand, there is a simple
one-to-one correspondence  between the set of  fixed points of $\Phi$ and the set of the reverse
plane partitions of shape $\delta_{n+2k}/\delta_n$. This   completes the proof  of   Conjecture \ref{conjecture}.

\section{Staircase alternating arrays}

In this section, we introduce the structure of staircase alternating arrays. We show that the generating function of staircase alternating arrays is equal to the determinant in
\eqref{mainresult}. We  define an operation on staircase alternating arrays.
 This operation will paly a key role in the construction of the involution $\Phi$ on staircase alternating arrays. We also give some properties on this operation.

A staircase alternating array is an array in which the entries in each row is an alternating word subject to
 certain conditions. As a generalization of alternating permutations, an alternating word is a word $W=a_1a_2\cdots a_m$ of nonnegative integers such that
\[a_1\geq a_2 \leq a_3\geq a_4\leq \cdots .\]
For a positive integer $k$, a $k$-staircase alternating array $\Pi$ is an array of nonnegative integers  such that
 \begin{itemize}
\item[(1)] For $1\leq i\leq k$, the entries in  $i$-th row form  an alternating word;

\item[(2)] For $1\leq i\leq k-1$, the last entry in row $i+1$ is two elements right to the last entry in row $i$.
\end{itemize}
For example, below is a $3$-staircase alternating array:
\begin{equation}\label{1}
\begin{array}{cccccccccccccc}
& &3&3&4&2&2\\
5&4&6&5&7&4&5&4&5\\
&&&&1&1&3&1&1&1&3.
\end{array}
\end{equation}
If there is no confusion occurring, we also represent  $\Pi$  by a $k$-tuple $(W_1,W_2,\ldots,W_k)$ of alternating words where $W_i$
is the alternating word in the $i$-th row of $\Pi$.

Let $\Pi=(W_1,W_2,\ldots,W_k)$ be a $k$-staircase alternating array. For $1\leq i\leq k$, assume that the $i$-th row
in $\Pi$ has $\ell_i$ elements. Let $n$ be a positive integer. We say that $\Pi$ is of order $n$ if
\begin{equation}\label{XXX}
\{\ell_i-2i+2\,|\, 1\leq i\leq k\}=\{2n+1,2n+3,\ldots,2n+2k-1\}.
\end{equation}
Intuitively, if we draw  a vertical line right after the last element in the first row of $\Pi$, then for each $1\leq i\leq k$, the value $\ell_i-2i+2$ in the set  \eqref{XXX} is the number
of entries in row $i$ that lie to the left of this line.

For example, for the  3-staircase alternating array $\Pi$ in   \eqref{1}, we see that $\ell_1=5, \ell_2=9,$ and $\ell_3=7$ and thus
$\{\ell_1,\ell_2-2,\ell_3-4\}=\{3,5,7\}$. Hence the order of $\Pi$ is 1.

In view of \eqref{XXX}, we can associate each $k$-staircase alternating array of order $n$ with a unique  permutation
$\sigma=\sigma_1\sigma_2\cdots \sigma_k$ on $\{1,2,\ldots,k\}$ by setting
\begin{equation}\label{Per}
\sigma_i=\frac{\ell_i'-2n+1}{2},
\end{equation}
where $\ell_i'=\ell_i-2i+2$.
For example, the    permutation associated to the array in  \eqref{1}  is
$\sigma=231$.
Note  that
for $1\leq i< j\leq k$,
$\sigma_i>\sigma_j$   if and only if
the first entry in  row $i$ of $\Pi$ is  to the left of the first entry in  row  $j$ of $\Pi$.
We also define  the sign $\mathrm{sgn}(\Pi)$ of $\Pi$ to be the sign of
the  associated  permutation $\sigma$, namely,
\[\mathrm{sgn}(\Pi)=(-1)^{\mathrm{inv}(\sigma)},\]
where
\[\mathrm{inv}(\sigma)=|\{(i,j)\,|\,1 \leq i<j\leq k, \sigma_i>\sigma_j\}|\]
 is the number of inversions of $\sigma$.

Let $A(n,k)$ be the set of $k$-staircase alternating arrays of order $n$.
For $\Pi\in A(n,k)$, let $|\Pi|$ represent  the weight of $\Pi$, that is, the sum of entries in $\Pi$.
The following theorem shows that the (signed) weight generating function of  staircase alternating arrays in $A(n,k)$ equals the  determinant in
 \eqref{mainresult}.

\begin{thm}\label{right}
For any positive integers $n$ and $k$,
\begin{equation}
\sum_{\Pi\in A(n,k)} {\mathrm{sgn}(\Pi)}q^{|\Pi|}=\det\left[\widetilde{E}^{*}_{2(n+i+j)-3}(q)
\right]_{i,j=1}^{k}.
\end{equation}
\end{thm}

\pf
For any fixed permutation $\sigma=\sigma_1\sigma_2\cdots \sigma_k$ on $\{1,2,\ldots,k\}$, let $A_\sigma$ denote  the subset of $A(n,k)$   consisting of arrays in  $A(n,k)$ to which the associated permutation is $\sigma$. Equivalently, a  $k$-tuple   $(W_1,W_2,\ldots,W_k)$ of alternating words
belongs to $A_\sigma$ if and only if the length of  $W_i$ is $2(n+i+\sigma_i)-3$.
Clearly,  $A(n,k)$ is a disjoint union of $A_\sigma$. Hence,
\begin{equation}\label{trans}
\sum_{\Pi\in A(n,k)} {\mathrm{sgn}(\Pi)}q^{|\Pi|}=\sum_{\sigma\in S_k}{\mathrm{sgn}(\sigma)}\sum_{\Pi\in A_\sigma}q^{|\Pi|},
\end{equation}
where $S_k$ is the set of
permutations on $\{1,2,\ldots,k\}$.

\begin{figure}[h t]
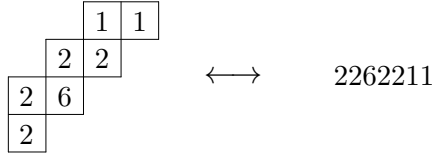

\centertexdraw{ \drawdim mm \linewd 0.1 \setgray 0

\move(90 5)\lvec(95 5) \move(90 10)\lvec(100 10)
\move(90 15)\lvec(105 15) \move(95 20)\lvec(110 20)
\move(100 25)\lvec(110 25)

\move(90 5)\lvec(90 15)\move(95 5)\lvec(95 20)
\move(100 10)\lvec(100 25)\move(105 15)\lvec(105 25)
 \move(110 20)\lvec(110 25)

\textref h:C v:C \htext(92.5 7.4){{\small 2}}

\textref h:C v:C \htext(92.5 12.4){{\small 2}}
\textref h:C v:C \htext(97.5 12.4){{\small 6}}

\textref h:C v:C \htext(97.5 17.4){{\small 2}}
\textref h:C v:C \htext(102.5 17.4){{\small 2}}

\textref h:C v:C \htext(102.5 22.4){{\small 1}}
\textref h:C v:C \htext(107.5 22.4){{\small 1}}

\textref h:C v:C \htext(120 15){$\longleftrightarrow$}
\textref h:C v:C \htext(140 15.2){{\small 2262211}}

} \caption{A reverse plane partition in the ribbon  $\delta_{n+2}/\delta_n$ and an alternating word}\label{corres}
\end{figure}

We   claim that for any permutation $\sigma=\sigma_1\sigma_2\cdots \sigma_k\in S_k$,
\begin{equation}\label{naog}
\sum_{\Pi\in A_\sigma}q^{|\Pi|}=\prod_{i=1}^k\widetilde{E}^{*}_{2(n+i+\sigma_i)-3}(q).
\end{equation}
The above claim can be shown as follows.
As mentioned in Introduction, Morales, Pak and Panova \cite{MPP} derived  that
\begin{equation}\label{ribbon}
\sum_{\pi\in \mathrm{RPP}(\delta_{n+2}/\delta_n)}{q^{|\pi|}}=\widetilde{E}^{*}_{2n+1}(q).
\end{equation}
On the other hand,  there is a one-to-one correspondence between
the set of reverse plane partitions of shape $\delta_{n+2}/\delta_n$
and the set of alternating words of length $2n+1$.
This is simply achieved by reading the entries of a reverse plane partition  in the ribbon  $\delta_{n+2}/\delta_n$
from the bottom left square to the top right square, see Figure \ref{corres}
for an example.
Thus, \eqref{ribbon} can be rewritten as
\begin{equation}\label{alter}
\sum_{W}{q^{|W|}}=\widetilde{E}^{*}_{2n+1}(q),
\end{equation}
where the sum is over alternating words of length $2n+1$.
Notice that
\begin{equation}\label{PP}
\sum_{\Pi\in A_\sigma}q^{|\Pi|}=\prod_{i=1}^k \sum_{W_i} q^{|W_i|},
\end{equation}
where $W_i$   runs over alternating words of length $2(n+i+\sigma_i)-3$.
Together with \eqref{alter},  we reach  the assertion  in \eqref{naog}.

Combining \eqref{trans} and \eqref{naog}, we obtain that
\begin{align*}
\sum_{\Pi\in A(n,k)} {\mathrm{sgn}(\Pi)}q^{|\Pi|}&=
\sum_{\sigma\in S_k}\mathrm{sgn}(\sigma)\sum_{\Pi\in A_\sigma}
q^{|\Pi|}\\[5pt]
&=\sum_{\sigma\in S_k} \mathrm{sgn}(\sigma) \prod_{i=1}^n\widetilde{E}^{*}_{2(n+i+\sigma_i)-3}(q)\\[5pt]
&=\det\left[\widetilde{E}^{*}_{2(n+i+j)-3}(q)
\right]_{i,j=1}^{k},
\end{align*}
as desired.
\qed

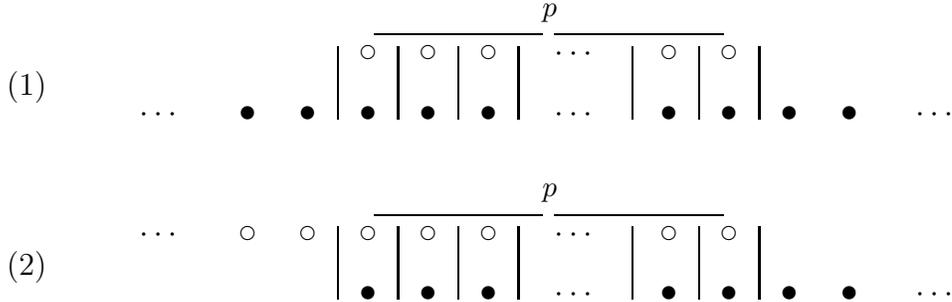
\begin{figure}[h t]
\setlength{\unitlength}{0.8mm}
\begin{center}
\begin{picture}(50,50)
\put(-20,10){\circle{2}}\put(-10,10){\circle{2}}
\put(0,10){\circle{2}} \put(10,10){\circle{2}}\put(20,10){\circle{2}}\put(31,8.5){$\cdots$}
\put(50,10){\circle{2}}\put(60,10){\circle{2}}\put(-38,8.5){$\cdots$}

\put(0,0){\circle*{2}} \put(10,0){\circle*{2}}\put(20,0){\circle*{2}}\put(31,-1.5){$\cdots$}
\put(50,0){\circle*{2}}\put(60,0){\circle*{2}}\put(70,0){\circle*{2}}
\put(80,0){\circle*{2}}\put(91,-1.5){$\cdots$}

\put(-5,-1){\line(0,1){12}}
\put(5,-1){\line(0,1){12}}\put(15,-1){\line(0,1){12}}\put(25,-1){\line(0,1){12}}
\put(44,-1){\line(0,1){12}}\put(55,-1){\line(0,1){12}}\put(65,-1){\line(0,1){12}}

\put(0,12){\line(1,1){1}}\put(59,13){\line(1,-1){1}}
\put(29,13){\line(1,1){1}}\put(30,14){\line(1,-1){1}}
\put(1,13){\line(1,0){28}}\put(31,13){\line(1,0){28}}
\put(29,16){{\small $p$}}

\put(-20,30){\circle*{2}}\put(-10,30){\circle*{2}}
\put(0,40){\circle{2}} \put(10,40){\circle{2}}\put(20,40){\circle{2}}\put(31,38.5){$\cdots$}
\put(50,40){\circle{2}}\put(60,40){\circle{2}}\put(-38,28.5){$\cdots$}

\put(0,30){\circle*{2}} \put(10,30){\circle*{2}}\put(20,30){\circle*{2}}\put(31,28.5){$\cdots$}
\put(50,30){\circle*{2}}\put(60,30){\circle*{2}}\put(70,30){\circle*{2}}
\put(80,30){\circle*{2}}\put(91,28.5){$\cdots$}

\put(-5,29){\line(0,1){12}}
\put(5,29){\line(0,1){12}}\put(15,29){\line(0,1){12}}\put(25,29){\line(0,1){12}}
\put(44,29){\line(0,1){12}}\put(55,29){\line(0,1){12}}\put(65,29){\line(0,1){12}}

\put(0,42){\line(1,1){1}}\put(59,43){\line(1,-1){1}}
\put(29,43){\line(1,1){1}}\put(30,44){\line(1,-1){1}}
\put(1,43){\line(1,0){28}}\put(31,43){\line(1,0){28}}
\put(29,46){{\small $p$}}

\put(-60,3){(2)}\put(-60,33){(1)}
\end{picture}
\end{center}
\caption{An illustration of positions  for the  case of $\sigma_i>\sigma_j$.}\label{f-1}
\end{figure}

We next define an operation on
staircase alternating arrays in $A(n,k)$.
For $\Pi\in A(n,k)$, let $\sigma=\sigma_1\sigma_2\cdots \sigma_k$ be  the permutation associated  to $\Pi$.   Let $W_i$ ($1\leq i\leq k$) denote   the $i$-th row of $\Pi$. For $1\leq i<j\leq k$, consider the elements in   $W_i$ and  $W_j$ that  overlap.
Assume  that there are $p$ such elements in each row. It follows from the definition of a $k$-staircase alternating array of order $n$ that $p$
must be an odd number. These elements
give rise to $p+1$ positions in row $i$ and row $j$. As an illustration, we signify these positions in Figure \ref{f-1} by short vertical lines, where  empty circles  represent the elements in row $i$ and
solid circles represent the elements in row $j$, and  (1) and (2) respectively correspond to  the cases  $\sigma_i<\sigma_j$ and $\sigma_i>\sigma_j$.

For each of the $p+1$ positions, we can get two new rows $W_i'$ and $W_j'$ by exchanging the elements of
$W_i$ and $W_j$  before this position. If  both $W_i'$ and $W_j'$ are still alternating words, then  the position  is called a  cutting  position. We say that $W_i$ and $W_j$
are  transposable if they have a cutting position. If  $W_i$ and $W_j$
are  transposable, then we can define a new
$k$-staircase alternating array, denoted  $T_{i,j}(\Pi)$, in $A(n,k)$ as follows.
Locate  the leftmost cutting position  of $W_i$ and $W_j$. Define  $T_{i,j}(\Pi)$ to be  the $k$-staircase alternating array obtained from
$\Pi$ by exchanging the elements of $W_i$ and $W_j$ before this cutting position.
By the above construction, it is easily seen that
\[T_{i,j}(T_{i,j}(\Pi))=\Pi.\]

Let us illustrate the operator $T_{i,j}$ by an  example. Let $\Pi$ be the staircase
alternating array in \eqref{1}. Consider the first row and the last row of $\Pi$.
There are two cutting positions  as illustrated in
Figure \ref{ff3}.
\begin{figure}[h t]
\setlength{\unitlength}{0.8mm}
\begin{center}
\begin{picture}(45,12)

\put(-20,5){3}\put(-10,5){3}
\put(0,5){4} \put(10,5){2}\put(20,5){2}

\put(0,-5){1} \put(10,-5){1}\put(20,-5){3}
\put(30,-5){1} \put(40,-5){1}\put(50,-5){1} \put(60,-5){3}

\put(16,-5){\line(0,1){14}}\put(26,-5){\line(0,1){14}}

\end{picture}
\end{center}
\caption{The cutting positions of two rows}\label{ff3}
\end{figure}
Thus, $T_{1,3}(\Pi)$ is the array obtained from
$\Pi$ by exchanging the elements of row 1 and row 3 before the left cutting position:
\begin{equation*}
\begin{array}{cccccccccccccc}
& &&&1&1&2\\[5pt]
5&4&6&5&7&4&5&4&5\\[5pt]
&&3&3&4&2&3&1&1&1&3
\end{array}
\end{equation*}

The rest of this section is denoted to several lemmas on the properties of the operation defined above.
Let $\Pi\in A(n,k)$, and $\sigma=\sigma_1\sigma_2\cdots \sigma_k$ be the  permutation associated to $\Pi$. Assume that  row $i$ and row $j$ of $\Pi$ are transposable, where  $1\leq i<j\leq k$.  We first give characterizations of   the first cutting position of row $i$ and row $j$. We  need to distinguish  the case $\sigma_i<\sigma_j$ and the case $\sigma_i>\sigma_j$.

We first consider the case $\sigma_i<\sigma_j$. In this case, row $i$ and row $j$ are of the following form:
\begin{equation}\label{z-1}
\begin{array}{cccccccccccccccccc}
& & &a_1&\ldots&a_{s}&\ldots &a_{2m+1}\\[5pt]
 \ldots&  b_1&b_2&b_3&\ldots&b_{s+2}&\ldots&b_{2m+3}&b_{2m+4}&b_{2m+5}&\ldots.
\end{array}
\end{equation}
As shown in Figure \ref{cutting-posi}, a cutting position is either right before or right after  an element    $a_{2t+1}$ ($0\leq t\leq m$) at  odd position of row $i$.
\begin{figure}[h t]
\setlength{\unitlength}{0.35mm}
\begin{center}
\begin{picture}(350,100)
%

\put(0,10){$
\begin{array}{cccccccccccccccccc}
& & &a_1&\cdots&a_{2t}&a_{2t+1}&a_{2t+2}&\cdots &a_{2m+1}\\[5pt]
 \cdots&  b_1&b_2&b_3&\cdots&b_{2t+2}&b_{2t+3}&b_{2t+4}&\cdots&b_{2m+3}&b_{2m+4}&b_{2m+5}&\cdots
\end{array}
$}
\put(-30,10){(2)}

\put(0,70){$
\begin{array}{cccccccccccccccccc}
& & &a_1&\cdots&a_{2t}&a_{2t+1}&a_{2t+2}&\cdots &a_{2m+1}\\[5pt]
 \cdots&  b_1&b_2&b_3&\cdots&b_{2t+2}&b_{2t+3}&b_{2t+4}&\cdots&b_{2m+3}&b_{2m+4}&b_{2m+5}&\cdots
\end{array}
$}
\put(-30,70){(1)}

\put(181,-3){\line(0,1){28}}\put(147,57){\line(0,1){28}}
\end{picture}
\end{center}\caption{A cutting position of row $i$ and row $j$ in the case $\sigma_i<\sigma_j$.}
\label{cutting-posi}
\end{figure}

\begin{lem}\label{bu-1}
Consider a cutting position of row $i$ and row $j$ as illustrated in Figure
\ref{cutting-posi}.
\begin{itemize}
\item[(1)] A cutting position of row $i$ and row $j$ right before     $a_{2t+1}$
is the first cutting position if and only if for $0\leq s\leq t-1$,
\begin{equation}\label{pp-2}
b_{2s+2}>a_{2s+1}<b_{2s+4}
\end{equation}
and
\begin{equation}\label{pp-1}
a_{2t+1}\geq b_{2t+2}.
\end{equation}

\item[(2)] A cutting position of row $i$ and row $j$ right after    $a_{2t+1}$
is the first cutting position if and only if for $0\leq s\leq t-1$,
\begin{equation}
b_{2s+2}>a_{2s+1}<b_{2s+4}
\end{equation}
and
\begin{equation}
b_{2t+2}>a_{2t+1}\geq b_{2t+4}.
\end{equation}
\end{itemize}
\end{lem}

\pf We shall only give a proof for the  assertion when the first cutting position is right before $a_{2t+1}$ as illustrated in (1) of Figure \ref{cutting-posi}. The argument for the other assertion  can be carried out in a similar manner.
When $t=0$, the assertion is obvious. We now consider the case for $t>0$.
Since the position right before $a_{2t+1}$ is a cutting position, the word
\[\begin{array}{cccccccccccccccccc}
 b_1&b_2&b_3&\cdots&b_{2t+2}&a_{2t+1}&a_{2t+2}&\cdots &a_{2m+1}
\end{array}\]
is alternating, and thus   \eqref{pp-1} holds. It remains to verify \eqref{pp-2}.
We only  check the  case for  $s=0$, namely,
\begin{equation}\label{ti-1}
b_2>a_1<b_4.
\end{equation}
The same argument   applies to the case for
 $s=1,2,\ldots, t-1$.

Suppose to the contrary that
$ a_1\geq b_2.$
Then it is easily seen that  the position right before $a_1$ is a cutting position, leading to a contradiction. So we get $a_1>b_2$.
We proceed to show that
$a_1<b_4$. Suppose otherwise that
\begin{equation}\label{ti-2}
a_1\geq b_4.
\end{equation}
Since the words in  row $i$ and row $j$ are alternating, we see that
 $b_3\geq b_2$ and $a_1\geq a_2$, which, together  with that face that $b_2>a_1$, leads to
 \begin{equation}\label{ti-3}
 b_3>a_2.
 \end{equation}
 By \eqref{ti-2} and \eqref{ti-3}, we see that
 the position right after $a_1$ is a cutting position,
which also leads to a contradiction. So we  have
 $a_1<b_4$. This verifies  \eqref{ti-1}.

We now prove  the reverse direction, that is,
if row $i$ and row $j$ satisfy the relations  in  \eqref{pp-2} and \eqref{pp-1}, then the position right before $a_{2t+1}$ is the first cutting position of row $i$ and row $j$.
 By \eqref{pp-2}, it follows that for any $0\leq s\leq t-1$, neither the  position
  right before $a_{2s+1}$ nor the position right after $a_{2s+1}$ is a cutting position. It remains to verify that the position
 right before  $a_{2t+1}$ is a cutting position. In view of \eqref{pp-1},  it suffices to show that
\begin{equation}\label{ti-4}
a_{2t}\leq b_{2t+3}.
\end{equation}
Since the words in row $i$ and row $j$ are alternating, it follows that
\begin{equation}\label{ti-5}
a_{2t}\leq a_{2t-1}\  \ \text{and}\ \ b_{2t+2}\leq b_{2t+3}.
\end{equation}
On the other hand, by \eqref{pp-2} we
see that
\begin{equation}\label{ti-6}
a_{2t-1}< b_{2t+2}.
\end{equation}
Combining \eqref{ti-5} and \eqref{ti-6}, we get \eqref{ti-4}. This completes the proof.
\qed

We next consider the  case $\sigma_i>\sigma_j$. In this case, row $i$ and row $j$ look like
\begin{equation}\label{z-2}
\begin{array}{cccccccccccccccccc}
 \ldots& a_1& a_2&a_3&\ldots&a_{s+2}&\ldots &a_{2m+3}\\[5pt]
&&&b_1&\ldots&b_{s}&\ldots&b_{2m+1}&b_{2m+2}&b_{2m+3}&\ldots.
\end{array}
\end{equation}
Moreover, as shown in Figure \ref{cutting-posi-2}, a cutting position is either right before or right after  an element    $b_{2t+1}$ ($0\leq t\leq m$) at  odd position of row $j$.
\begin{figure}[h t]
\setlength{\unitlength}{0.35mm}
\begin{center}
\begin{picture}(350,100)
%

\put(0,10){$
\begin{array}{cccccccccccccccccc}
\ldots&a_1 & a_2&a_3&\ldots&a_{2t+2}&a_{2t+3}&a_{2t+4}&\ldots &a_{2m+3}\\[5pt]
&&&b_1&\ldots&b_{2t}&b_{2t+1}&b_{2t+2}&\ldots&b_{2m+1}&b_{2m+2}&b_{2m+3}&\ldots.
\end{array}
$}
\put(-30,10){(2)}

\put(0,70){$
\begin{array}{cccccccccccccccccc}
\ldots&a_1 & a_2&a_3&\ldots&a_{2t+2}&a_{2t+3}&a_{2t+4}&\ldots &a_{2m+3}\\[5pt]
&&&b_1&\ldots&b_{2t}&b_{2t+1}&b_{2t+2}&\ldots&b_{2m+1}&b_{2m+2}&b_{2m+3}&\ldots.
\end{array}
$}
\put(-30,70){(1)}

\put(180,-3){\line(0,1){28}}\put(145,57){\line(0,1){28}}
\end{picture}
\end{center}\caption{A cutting position of row $i$ and row $j$
in the case $\sigma_i>\sigma_j$.}
\label{cutting-posi-2}
\end{figure}
The following lemma gives a characterization of the first cutting position of row $i$ and row $j$ in the case $\sigma_i>\sigma_j$. The proof  is similar  to that in
Lemma \ref{bu-1} for the case $\sigma_i<\sigma_j$, and thus is omitted.

\begin{lem}\label{bu-2}
Consider a cutting position of row $i$ and row $j$ as illustrated in Figure
\ref{cutting-posi-2}.
\begin{itemize}
\item[(1)] A cutting position of row $i$ and row $j$ right before     $b_{2t+1}$
is the first cutting position if and only if for $0\leq s\leq t-1$,
\begin{equation}\label{qq-1}
a_{2s+2}>b_{2s+1}<a_{2s+4}
\end{equation}
and
\begin{equation}\label{qq-2}
b_{2t+1}\geq a_{2t+2}.
\end{equation}
\item[(2)] A cutting position of row $i$ and row $j$ right after     $b_{2t+1}$
  is the first cutting position if and only if for $0\leq s\leq t-1$,
\begin{equation}\label{qq-3}
a_{2s+2}>b_{2s+1}<a_{2s+4}
\end{equation}
and
\begin{equation}\label{qq-4}
a_{2t+2}>b_{2t+1}\geq a_{2t+4},
\end{equation}
where, in the case $t=m$, we set $a_{2m+4}=0$.
\end{itemize}
\end{lem}
%
%

Based on   Lemma \ref{bu-2}, we can prove the third lemma which asserts that
if $\sigma_i>\sigma_j$, then row $i$ and row $j$ must be transposable.

\begin{lem}\label{l1}
Let $\Pi\in A(n,k)$ and  $\sigma=\sigma_1\sigma_2\cdots \sigma_k$ be the permutation associated to $\Pi$. For $1\leq i<j\leq k$, if $\sigma_i>\sigma_j$, then row $i$ and row $j$ of $\Pi$ are transposable.
\end{lem}

\pf Assume that row $i$ and row $j$ are as illustrated below:
\begin{equation*}
\begin{array}{cccccccccccccccccc}
 \ldots& a_1& a_2&a_3&\ldots&a_{s+2}&\ldots &a_{2m+3}\\[5pt]
&&&b_1&\ldots&b_{s}&\ldots&b_{2m+1}&b_{2m+2}&b_{2m+3}&\ldots.
\end{array}
\end{equation*}
We claim  that there exists an index $t$ with $0\leq t\leq m $ such that
$b_{2t+1}$ satisfies either \eqref{qq-2} or \eqref{qq-4}.
If this is not the case, then for any   $0\leq t\leq m $,
\[a_{2t+2}>b_{2t+1}<a_{2t+4}.\]
In particular, we have $b_{2m+1}<a_{2m+4}=0$, leading to a contradiction.
This verifies the claim.

Let   $p$ be the smallest index such that
$b_{2p+1}$ satisfies either \eqref{qq-2} or \eqref{qq-4}. Then, for $0\leq t\leq p-1$,
\begin{equation}\label{ooq}
a_{2t+2}>b_{2t+1}< a_{2t+4}.
\end{equation}
If $b_{2p+1}$ satisfies  \eqref{qq-2}, by (1) of Lemma \ref{bu-2}  we see that the position right before $b_{2p+1}$ is a cutting position. Otherwise,
  $b_{2p+1}$ satisfies  \eqref{qq-4}. It follows from (2) of Lemma \ref{bu-2} that   the position right after $b_{2p+1}$ is a cutting position.
Hence, in both cases, row $i$ and row $j$ are transposable. This completes the proof.
\qed

Finally, combining Lemma \ref{bu-1} and Lemma \ref{l1}, we are led to
the following  characterization on when row $i$ and row $j$ are not transposable.

\begin{lem}\label{ll1}
Let $\Pi\in A(n,k)$ and  $\sigma=\sigma_1\sigma_2\cdots \sigma_k$ be the permutation associated to $\Pi$. Then, for $1\leq i<j\leq k$, row $i$ and row $j$ in $\Pi$ are not transposable if and only if
\begin{itemize}
\item[(1)] $\sigma_i<\sigma_j$;
\item[(2)] Assume that row $i$ and row $j$ are of the following form
\begin{equation*}
\begin{array}{cccccccccccccccccc}
& & &a_1&\ldots&a_{2t}&a_{2t+1}&a_{2t+2}&\ldots &a_{2m+1}\\[5pt]
 \ldots&  b_1&b_2&b_3&\ldots&b_{2t+2}&b_{2t+3}&b_{2t+4}&\ldots&b_{2m+3}&b_{2m+4}&b_{2m+5}&\ldots.
\end{array}
\end{equation*}
Then,
for $0\leq t\leq m$,
\begin{equation}\label{keyineq}
b_{2t+2}>a_{2t+1}<b_{2t+4}.
\end{equation}
\end{itemize}

\end{lem}

\section{Proof of  Conjecture \ref{conjecture}}

In this section, we provide a proof of Conjecture \ref{conjecture}. We construct  a weight preserving involution $\Phi$ on the set
 $A(n,k)$ of $k$-staircase
alternating arrays of order $n$. The fixed points of $\Phi$
are   arrays of $A(n,k)$ that have no transposable  rows.
By Lemma \ref{ll1}, the associated permutation to
each fixed point  of $\Phi$ is the identity permutation, and so  each  fixed point of $\Phi$ has a positive sign. On the other hand, every non-fixed point of $\Pi$ has an opposite sign with
$\Phi(\Pi)$. Hence we have
\begin{equation}\label{tran}
\sum_{\Pi\in A(n,k)} {\mathrm{sgn}(\Pi)}q^{|\Pi|}=\sum_{\Pi\in FA(n,k)} q^{|\Pi|},
\end{equation}
where $FA(n,k)$ is the set of fixed points of $\Phi$.
Moreover, as will be seen in Theorem \ref{La}, there is an obvious  bijection between the fixed points of
$\Phi$ and reverse plane partitions of shape $\delta_{n+2k}/\delta_n$, which implies the following relation
\begin{equation}\label{vt}
q^{N}\sum_{\pi\in \mathrm{RPP}(\delta_{n+2k}/\delta_n)}{q^{|\pi|}}=\sum_{\Pi\in FA(n,k)} q^{|\Pi|}.
\end{equation}
Combining \eqref{tran}, \eqref{vt} and Theorem \ref{right}, we are led to a proof of Conjecture \ref{conjecture}.

We first describe a weight preserving map $\Phi$ on  $A(n,k)$, and then show that $\Phi$ is an involution.

\noindent
\textbf{A map $\Phi$ on $A(n,k)$:}
Let $\Pi\in A(n,k)$, and let $\sigma=\sigma_1\sigma_2\cdots \sigma_k$ be the permutation associated to $\Pi$. If,
for any $1\leq m<k$, row $m$ and row $m+1$ of $\Pi$ are not transposable, then we set $\Pi$ to be a fixed point of $\Phi$, namely,
\[\Phi(\Pi)=\Pi.\]
Otherwise, there is an index $1\leq m<k$ such that row $m$ and row $m+1$ of $\Pi$ are transposable.
Assume that  $m_0$ is the smallest such index. Note  that by Lemma \ref{ll1},
\[\sigma_1<\sigma_2<\cdots<\sigma_{m_0}.\]
To define $\Phi(\Pi)$, we consider two cases.

\noindent
(A): $\sigma_{m_0}<\sigma_{m_0+1}$. Set
\[\Phi(\Pi)=T_{m_0,m_0+1}(\Pi).\]

\noindent
(B): $\sigma_{m_0}>\sigma_{m_0+1}$. Let $1\leq s\leq m_0$ be the smallest index such $\sigma_s>\sigma_{m_0+1}$. By Lemma \ref{l1}, row $s$ and row $m_0+1$ are transposable. To define
$\Phi(\Pi)$, consider row $s-1$ and row $s$ in $T_{s,m_0+1}(\Pi)$.
There are  two subcases.

\noindent
(I): Row $s-1$ and row $s$ of $T_{s,m_0+1}(\Pi)$ are transposable.
As will be seen in Lemma \ref{Lemma-2}, row $s-1$ and row $m_0+1$ in $\Pi$ are transposable. Set
\[\Phi(\Pi)=T_{s-1,m_0+1}(\Pi).\]

\noindent
(II): Row $s-1$ and row $s$ of $T_{s,m_0+1}$ are not transposable.
Set
\[\Phi(\Pi)=T_{s,m_0+1}(\Pi).\]

The main theorem in this section  asserts that $\Phi$ is an involution.

\begin{thm} \label{invo}
The map $\Phi$ is  a weight preserving involution on $A(n,k)$.
\end{thm}

Clearly, $\Phi$ preserves the weight. To prove Theorem \ref{invo}, we need to show that $\Phi^2(\Pi)=\Pi$ for any
$\Pi\in A(n,k)$. This is obvious  when  $\Pi$ is a fixed point. To conclude $\Phi^2(\Pi)=\Pi$  for a non-fixed point $\Pi$, we need  three lemmas.
This first lemma is used to verify  that if $\Pi$ satisfies the condition in (A), then
$\Phi^2(\Pi)=\Pi$.

\begin{lem}\label{lem1}
Let $\Pi\in A(n,k)$, and   $\sigma=\sigma_1\sigma_2\cdots \sigma_k$ be the permutation associated to $\Pi$. Let $1\leq i<j<\ell\leq k$ be indices such that $\sigma_i<\sigma_j<\sigma_\ell$.¡¡Assume that  row $i$ and row $j$  of $\Pi$ are not transposable, and that row $j$ and row $\ell$  of $\Pi$ are transposable. Then, row $i$ and row $j$ in $T_{j,\ell}(\Pi)$ are not transposable.
\end{lem}

\pf Assume that row $i$, row $j$ and row $\ell$ of $\Pi$ look like
\[
\begin{array}{cccccccccccccccccc}
&& & & &a_1&\cdots&a_{2t}&a_{2t+1}&\cdots&a_{2m+1}\\
&&&b_1&\cdots&b_{2s+1}&\cdots&b_{2(s+t)}&b_{2(s+t)+1}&\cdots&b_{2(m+s)+1}&b_{2(m+s)+2}&\cdots\\
\cdots&c_1&c_2&c_3&\cdots&c_{2s+3}&\cdots&c_{2(s+t)+2}&c_{2(s+t)+3}&\cdots&c_{2(m+s)+3}&c_{2(m+s)+4}&\cdots
\end{array}
\]
Since row $i$ and row $j$ are not transposable, it follows from Lemma \ref{ll1} that
\begin{equation}\label{aleq-1}
b_{2s}>a_1<b_{2s+2}> \cdots<b_{2(m+s)}>a_{2m+1}<b_{2(m+s)+2}.
\end{equation}
On the other hand, since row $j$ and row $\ell$ are transposable, according to the position of the first cutting position, we have two cases.

\noindent
Case 1:
The first cutting position is right before  $b_{2(s+t)+1}$. By Lemma \ref{bu-1}, we see that
\begin{equation}\label{aleq-2}
c_2>b_1<c_4>\cdots<c_{2(s+t)}>b_{2(s+t)-1}<c_{2(s+t)+2}.
\end{equation}
Thus, row $i$, row $j$ and row $\ell$ in $T_{j,\ell}(\Pi)$ are of the following form
\[
\begin{array}{cccccccccccccccccc}
&& & & &a_1&\cdots&a_{2t}&a_{2t+1}&\cdots&a_{2m+1}\\
\cdots&c_1&c_2&c_3&\cdots&c_{2s+3}&\cdots&c_{2(s+t)+2}&b_{2(s+t)+1}&\cdots&b_{2(m+s)+1}&b_{2(m+s)+2}&\cdots\\
&&&b_1&\cdots&b_{2s+1}&\cdots&b_{2(s+t)}&c_{2(s+t)+3}&\cdots&c_{2(m+s)+3}&c_{2(m+s)+4}&\cdots
\end{array}
\]
By \eqref{aleq-1} and \eqref{aleq-2}, it is easy to check that
\[
c_{2s+2}>a_1<c_{2s+4}>\cdots <c_{2(s+t)}>a_{2t-1}<c_{2(s+t)+2}\ \ \ \text{and}\ \ \ a_{2t+1}<c_{2(s+t)+2},
\]
which, together with \eqref{aleq-1}, implies that row $i$ and row $j$ of
$T_{j,\ell}(\Pi)$ satisfy the conditions  in Lemma \ref{ll1}. Hence
  row $i$  and row $j$ of $T_{j,\ell}(\Pi)$ are not transposable.

\noindent
Case 2:
The first cutting position  is right  after  $b_{2(s+t)+1}$.
Then  row $i$, row $j$ and row $\ell$ in $T_{j,\ell}(\Pi)$ are of the following form
\[
\begin{array}{cccccccccccccccccc}
&& & & &a_1&\cdots&a_{2t+1}&a_{2t+2}&\cdots&a_{2m+1}\\
\cdots&c_1&c_2&c_3&\cdots&c_{2s+3}&\cdots&c_{2(s+t)+3}&b_{2(s+t)+2}&\cdots&b_{2(m+s)+1}&b_{2(m+s)+2}&\cdots\\
&&&b_1&\cdots&b_{2s+1}&\cdots&b_{2(s+t)+1}&c_{2(s+t)+4}&\cdots&c_{2(m+s)+3}&c_{2(m+s)+4}&\cdots
\end{array}
\]
Similar to Case 1, it is easy to verify that  row $i$ and row $j$ of
$T_{j,\ell}(\Pi)$ satisfy the conditions  in Lemma \ref{ll1}, and thus
  row $i$  and row $j$ of  $T_{j,\ell}(\Pi)$ are not transposable.
  This completes the proof.
\qed

We now state the second lemma. This lemma ensures  that   $\Phi$ is well-defined if it satisfies the condition in  (I) of (B). On the other hand, we use this lemma to show that if  $\Pi$ satisfies the condition in (I) of (B), then
$\Phi^2(\Pi)=\Pi$.

\begin{lem}\label{Lemma-2}
Let $\Pi\in A(n,k)$, and let $\sigma=\sigma_1\sigma_2\cdots \sigma_k$ be the permutation associated to $\Pi$. Assume that $1\leq i<j<\ell\leq k$ are indices such that $\sigma_i<\sigma_\ell<\sigma_j$ and that row $i$  and row $j$  in $\Pi$ are not transposable. By Lemma \ref{l1}, row $j$ and row $\ell$ in $\Pi$
 are transposable. If row $i$  and row $j$ in $T_{j,\ell}(\Pi)$ are transposable,
 then  row $i$  and row $\ell$ in $\Pi$ are transposable. Moreover,
row $i$  and row $j$ in $T_{i,\ell}(\Pi)$  are not transposable.

\end{lem}

\pf By the assumption  that $\sigma_i<\sigma_\ell<\sigma_j$, row $i$, row $j$ and row $\ell$ of $\Pi$ are of the following form
\[
\begin{array}{cccccccccccccccccc}
&& & & &a_1&\cdots&a_{2t}&a_{2t+1}&\cdots&a_{2m+1}\\
\cdots&b_1&b_2&b_3&\cdots&b_{2s+3}&\cdots&b_{2(s+t)+2}
&b_{2(s+t)+3}&\cdots&b_{2(m+s)+3}&b_{2(m+s)+4}&\cdots\\
&&&c_1&\cdots&c_{2s+1}&\cdots&c_{2(s+t)}&c_{2(s+t)+1}
&\cdots&c_{2(m+s)+1}&c_{2(m+s)+2}&\cdots
\end{array}
\]
Since row $i$ and row $j$ of $\Pi$ are not transposable, by Lemma \ref{ll1} we see that
\begin{equation}\label{aleq-11}
b_{2s+2}>a_1<b_{2s+4}>\cdots<b_{2(m+s)+2}>a_{2m+1}<b_{2(m+s)+4}.
\end{equation}
On the other hand, since row $i$  and row $j$ in $T_{j,\ell}(\Pi)$ are transposable, the first cutting position of row $j$ and row $\ell$ is to the right of
  $c_{2s}$. There are  two cases to consider.

\noindent
Case 1: The first cutting position of row $i$ and row $j$ in $\Pi$  is  right before  $c_{2(s+t)+1}$. By Lemma \ref{bu-2}, we have
\begin{equation}\label{puke}
b_2>c_1<b_4>\cdots<b_{2(s+t)}>c_{2(s+t)-1}<b_{2(s+t)+2}.
\end{equation}
So, row $i$, row $j$ and row $\ell$ of  $T_{j,\ell}(\Pi)$ are of the following form
\[
\begin{array}{cccccccccccccccccc}
&& & & &a_1&\cdots&a_{2t}&a_{2t+1}&\cdots&a_{2m+1}\\
&&&c_1&\cdots&c_{2s+1}&\cdots&c_{2(s+t)}&b_{2(s+t)+3}&\cdots&b_{2(m+s)+3}&b_{2(m+s)+4}&\cdots\\
\cdots&b_1&b_2&b_3&\cdots&b_{2s+3}&\cdots&b_{2(s+t)+2}&c_{2(s+t)+1}&\cdots&c_{2(m+s)+1}&c_{2(m+s)+2}&\cdots.
\end{array}
\]
Keep in mind that row $i$ and row $j$ in $\Pi$ are not transposable, and that
 row $i$ and row $j$ in $T_{j,\ell}(\Pi)$ are   transposable. So  the first cutting position of row $i$ and row $j$ in $T_{j,\ell}(\Pi)$  is a position right before or right after $a_p$, where $p$ is an odd index  with $p\leq 2t+1$.  Clearly, this position is also the first cutting position  of  row $i$  and  row  $\ell$ of $\Pi$. Hence we conclude that row $i$ and row $\ell$ in
$\Pi$ are transposable.

We proceed to show that row $i$  and row $j$ in $T_{i,\ell}(\Pi)$ are not transposable. As we have discussed above, the   first cutting position  of  row $i$  and  row  $\ell$ of $\Pi$ is a position  right before or right after $a_p$, where $p$ is an odd index  with $p\leq 2t+1$. So row $i$ and row $j$ in $T_{i,\ell}(\Pi)$ have the following form
\[
\begin{array}{cccccccccccccccccc}
&&&c_1&\cdots&c_{2s+p-1}&a_p&\cdots&a_{2t+1}&\cdots&a_{2m+1}\\
\cdots&b_1&b_2&b_3&\cdots&b_{2s+p+1}&b_{2s+p+2}&\cdots
&b_{2(s+t)+3}&\cdots&b_{2(m+s)+3}&b_{2(m+s)+4}&\cdots
\end{array}
\]
By \eqref{aleq-11} and \eqref{puke}, it is easily seen that row $i$  and row $j$ in $T_{i,\ell}(\Pi)$ satisfy the conditions  in  Lemma \ref{ll1}, and thus they are not transposable.

\noindent
Case 2: The first cutting position of row $i$ and row $j$ in $\Pi$  is right after $c_{2(s+t)+1}$.
The argument for this case is analogous to that for Case 1, and thus is omitted.
This completes the proof.
\qed

Finally, we give the third  lemma. This lemma  allows us to deduce that if $\Pi$ satisfies the condition in (II) of Case (B), then
$\Phi^2(\Pi)=\Pi$.

\begin{lem}\label{Lemma-3}
Let $\Pi\in A(n,k)$, and let $\sigma=\sigma_1\sigma_2\cdots \sigma_k$ be the permutation associated to $\Pi$. Assume that $1\leq i<j<\ell\leq k$ are indices such that $\sigma_\ell<\sigma_i<\sigma_j$ and that
row $i$  and row $j$ in $\Pi$  are not transposable. By Lemma \ref{l1}, row $i$  and row $\ell$ in $\Pi$  are transposable, and row  $j$  and row $\ell$ in $T_{i,\ell}(\Pi)$ are transposable.
Then,
\begin{itemize}
 \item[(1)] Row $i$  and row $j$ in $T_{i,\ell}(\Pi)$ are not transposable;

 \item[(2)] Row $i$ and row $j$ in $T_{j,\ell}(T_{i,\ell}(\Pi))$ are transposable.
\end{itemize}
\end{lem}

\pf By the assumption that $\sigma_\ell<\sigma_i<\sigma_j$,  row $i$, row $j$ and row $\ell$ of $\Pi$ are of the following form
\[
\begin{array}{cccccccccccccccccc}
&&&a_1&\cdots&a_{2s+1}&\cdots&a_{2(s+t)}&a_{2(s+t)+1}&\cdots&a_{2m+1}
\\
\cdots&b_1&b_2&b_3&\cdots&b_{2s+3}&\cdots&b_{2(s+t)+2}&b_{2(s+t)+3}&\cdots&b_{2m+3}&b_{2m+4}&\cdots\\
&& & & &c_1&\cdots&c_{2t}&c_{2t+1}&\cdots&c_{2(m-s)+1}&c_{2(m-s)+2}&\cdots
\end{array}
\]
Since row $i$ and row $j$ of $\Pi$ are not transposable, it follows from Lemma \ref{ll1} that
\begin{equation}\label{la-1}
b_{2}>a_1<b_4>\cdots<b_{2m+2}>a_{2m+1}<b_{2m+4}.
\end{equation}
By Lemma \ref{l1}, row $j$ and row $\ell$ of $\Pi$ are transposable. According to the first cutting position of row $j$ and row $\ell$ in $\Pi$,
we consider the following  two cases.

\noindent
Case 1: The first cutting position of row $j$ and row $\ell$ in $\Pi$  is right before $c_{2t+1}$.  By Lemma \ref{bu-2}, we obtain that
\begin{equation}\label{la-2}
b_{2s+2}>c_1<b_{2s+4}>\cdots<b_{2(s+t)}>c_{2t-1}<b_{2(s+t)+2}
\end{equation}
and
\begin{equation}\label{La2}
c_{2t+1}\geq b_{2(s+t)+2}.
\end{equation}
By  \eqref{la-1} and \eqref{La2}, we get
\begin{equation}\label{la-3}
 c_{2t+1}\geq b_{2(s+t)+2} > a_{2(s+t)+1}\geq a_{2(s+t)}.
\end{equation}
In view of Lemma \ref{bu-2}, \eqref{la-3} implies that the first cutting position of row $i$ and row $\ell$ in $\Pi$ is a position right before $c_{i+1}$ where $i\leq 2t$.
Hence row $i$, row $j$ and row $\ell$ in $T_{i,\ell}(\Pi)$  look like
\[
\begin{array}{cccccccccccccccccc}
&& & & &c_1&\cdots&c_{i}&a_{2s+i+1}&\cdots&a_{2(s+t)+1}&\cdots&a_{2m+1}
\\
\cdots&b_1&b_2&b_3&\cdots&b_{2s+3}&\cdots&b_{2s+i+2}&b_{2s+i+3}&\cdots&b_{2(s+t)+3}&\cdots&b_{2m+3}&\cdots\\
&&&a_1&\cdots&a_{2s+1}&\cdots&a_{2s+i}&c_{i+1}&\cdots&c_{2t+1}&\cdots&c_{2(m-s)+1}&\cdots
\end{array}
\]
By \eqref{la-1} and \eqref{la-2}, it follows that row $i$ and row $j$
in $T_{i,\ell}(\Pi)$ satisfy the conditions in Lemma \ref{ll1}, and thus they are not transposable.

We next show that row $i$ and row $j$ in $T_{j,\ell}(T_{i,\ell}(\Pi))$ are transposable. Notice that the first cutting position of row $j$ and row $\ell$
in $T_{i,\ell}(\Pi)$ is the position right before $c_{2t+1}$. Hence row $i$, row $j$ and row $\ell$ in $T_{j,\ell}(T_{i,\ell}(\Pi))$ are of the following form
\[
\begin{array}{cccccccccccccccccc}
&& & & &c_1&\cdots&c_{i}&a_{2s+i+1}&\cdots&a_{2(s+t)+1}&\cdots&a_{2m+1}
\\&&&a_1&\cdots&a_{2s+1}&\cdots&a_{2s+i}&c_{i+1}&\cdots&b_{2(s+t)+3}
&\cdots&b_{2m+3}&\cdots\\
\cdots&b_1&b_2&b_3&\cdots&b_{2s+3}&\cdots&b_{2s+i+2}&b_{2s+i+3}&\cdots&c_{2t+1}&\cdots&c_{2(m-s)+1}&\cdots
\end{array}
\]
If $i$ is odd, then we see that $c_i>c_{i+1}$. Otherwise, if $i$ is even, then we
have $a_{2s+i+1}>a_{2s+i}$. In both cases, row $i$ and row $j$
in $T_{j,\ell}(T_{i,\ell}(\Pi))$  do not satisfy Condition (2) in  Lemma \ref{ll1}
and thus they are transposable.

\noindent
Case 2: The first cutting position of row $j$ and row $\ell$ in $\Pi$  is right after $c_{2t+1}$. The argument is analogous to that for Case 1, and is omitted here.
This completes the proof.

\qed

With the above three lemmas, we can now give a proof of Theorem \ref{invo}.

\noindent
{\it Proof of Theorem \ref{invo}.} Let $\Pi\in A(n,k)$, and $\sigma=\sigma_1\sigma_1\cdots \sigma_k$ be the associated permutation. If $\Pi$ is a fixed point of $\Phi$, then it is clear that $\Phi^(\Pi)=\Pi$. We next consider the case when $\Pi$ is not a fixed point.
Assume that  $m_0$ ($1\leq m_0<k$) is the smallest index such that row $m_0$ and row $m_0+1$ in $\Pi$ are transposable. By Lemma \ref{ll1},
\[\sigma_1<\sigma_2<\cdots<\sigma_{m_0}.\]
According to the construction of $\Phi$, we consider the following two cases.

\noindent
Case 1: $\sigma_{m_0}<\sigma_{m_0+1}$. In this case,
\[\Phi(\Pi)=T_{m_0,m_0+1}(\Pi).\]
By Lemma \ref{lem1}, we see that  $m_0$   is still the smallest index such that row $m_0$ and row $m_0+1$ in $\Phi(\Pi)$ are transposable. Moreover, $\Phi(\Pi)$ satisfies the condition (II) in (B). So we reach that
\[\Phi^2(\Pi)=\Phi(T_{m_0,m_0+1}(\Pi))=T_{m_0,m_0+1}(T_{m_0,m_0+1}(\Pi))=\Pi.\]

\noindent
Case 2: $\sigma_{m_0}>\sigma_{m_0+1}$. Let $1\leq s\leq m_0$ be the smallest index such $\sigma_s>\sigma_{m_0+1}$. There are    two subcases.

\noindent
(1): Row $s-1$ and row $s$ of $T_{s,m_0+1}$ are transposable.
In this case,
\[\Phi(\Pi)=T_{s-1,m_0+1}(\Pi).\]
By Lemma \ref{lem1}, it follows  that row $s-2$ and row $s-1$ in $\Phi(\Pi)$ are not transposable. Moreover, by Lemma \ref{Lemma-2} we see that row $s-1$ and row $s$ in $\Phi(\Pi)$ are not transposable. This implies that $m_0$   is also the smallest index such that row $m_0$ and row $m_0+1$ in $\Phi(\Pi)$ are transposable.
It is easily checked that  $\Phi(\Pi)$ satisfies the condition (II) in (B).
Hence we obtain that
\[\Phi^2(\Pi)=\Phi(T_{s-1,m_0+1}(\Pi))=T_{s-1,m_0+1}(T_{s-1,m_0+1}(\Pi))=\Pi.\]

\noindent
(2): Row $s-1$ and row $s$ of $T_{s,m_0+1}$ are not transposable.
In this case,
\[\Phi(\Pi)=T_{s,m_0+1}(\Pi).\]
By the assertion  (1) in Lemma \ref{Lemma-3},    row $s$ and row $s+1$ in $\Phi(\Pi)$ are not transposable. So, we see that  $m_0$   is  the smallest index such that row $m_0$ and row $m_0+1$ in $\Phi(\Pi)$ are transposable. By the assertion (2)
in Lemma \ref{Lemma-3}, it is easy to check  that
$\Phi(\Pi)$ satisfies the condition (I) in (B).
Hence we have
\[\Phi^2(\Pi)=\Phi(T_{s,m_0+1}(\Pi))=T_{s,m_0+1}(T_{s,m_0+1}(\Pi))=\Pi.\]
This completes the proof.
\qed

Notice that  every non-fixed point $\Pi$ has an opposite sign with $\Phi(\Pi)$.
Therefore, as a direct consequence of  Theorem \ref{invo}, we obtain that
\begin{equation}\label{Re}
\sum_{\Pi\in A(n,k)} {\mathrm{sgn}(\Pi)}q^{|\Pi|}=\sum_{\Pi\in FA(n,k)} {\mathrm{sgn}(\Pi)}q^{|\Pi|}.
\end{equation}
For a fixed point  $\Pi\in FA(n,k)$ of $\Phi$, let $\sigma=\sigma_1\sigma_2\cdots\sigma_k$ be the permutation
associated to $\Pi$.
 According to Lemma \ref{l1}, we see that $\sigma_1<\sigma_2<\cdots<\sigma_k$,
which implies that $\sigma=12\cdots k$ is the identity permutation.   Hence, \eqref{Re} can be rewritten as
\begin{equation}\label{RRe}
\sum_{\Pi\in A(n,k)} {\mathrm{sgn}(\Pi)}q^{|\Pi|}=\sum_{\Pi\in FA(n,k)} q^{|\Pi|}.
\end{equation}
The following theorem shows   that the generating function on the right-hand side of \eqref{RRe} equals the left-hand side of \eqref{mainresult}.

\begin{thm}\label{La}
For any positive integers $n$ and $k$,
\[\sum_{\Pi\in FA(n,k)}q^{|\Pi|}=q^{N}\sum_{\pi\in \mathrm{RPP}(\delta_{n+2k}/\delta_n)}{q^{|\pi|}},\]
where
\[N=\frac{k(k-1)(6n+8k-1)}{6},\]
\end{thm}

\pf We shall construct a bijection $\varphi$ from the set $FA(n,k)$ to the set
$\mathrm{RPP}(\delta_{n+2k}/\delta_n)$ such that for each $\Pi\in FA(n,k)$,
\begin{equation}\label{fina}
|\Pi|=|\varphi(\Pi)|+N.
\end{equation}

Let $\Pi\in FA(n,k)$ be a fixed point of $\Phi$, and let $W_i$ ($1\leq i\leq k$) denote the $i$-th row of $\Pi$. Since  the permutation associated to $\Pi$ is the identity permutation, by \eqref{Per} we see that $W_i$ has $2n+4i-3$ elements. More precisely, for $1\leq i\leq k-1$,   $W_i$ and  $W_{i+1}$    are of the following form
\begin{equation*}
\begin{array}{cccccccccccccc}
 & &a_1&a_2&a_3&\cdots&a_{2n+4i-3}\\[5pt]
b_1&b_{2}&b_{3}&b_{4}&b_{5}&\cdots&b_{2n+4i-1}&b_{2n+4i}&b_{2n+4i+1}
\end{array}
\end{equation*}
Since row $i$ and row $i+1$ in $\Pi$ are not transposable, by Lemma \ref{ll1} we see that
\begin{equation}\label{fi}
b_{2}>a_1<b_{4}>\cdots<b_{2n+4i-2}>a_{2n+4i-3}<b_{2n+4i}.
\end{equation}
This implies  that for $1\leq i\leq k$,
each element in  $W_i$ is greater than or equal to $i-1$.
Let $\Pi'$ be the $k$-staircase alternating array such that for
$1\leq i\leq k$,  the $i$-th row $W_i'$ of $\Pi'$ is    obtained from $W_i$ by subtracting each element by $i-1$.
Clearly,  $\Pi'$ corresponds to a reverse  plane
partition $\pi$ of shape $\delta_{n+2k}/\delta_{n}$ by putting $W_i'$ into the border strip  $\delta_{n+2i}/\delta_{n+2i-2}$ along the bottom left square to the top right square. Define
\[\varphi(\Pi)=\pi.\]
It is straightforward  to construct the reverse procedure of $\varphi$, and thus $\varphi$
is a bijection.

Notice that
\begin{align*}
|\Pi|&=|\Pi'|+\sum_{i=1}^k (i-1)(2n+4i-3)\\[5pt]
&=|\Pi'|+\sum_{i=0}^{k-1} i(2n+4i+1)\\[5pt]
&=|\Pi'|+(2n+1)\sum_{i=0}^{k-1} i(2n+1)+4\sum_{i=0}^{k-1} i^2\\[5pt]
&=|\Pi'|+(2n+1)\frac{ (k-1)k}{2}+4\frac{(k-1)k(2k-1)}{6}\\[5pt]
&=|\Pi'|+\frac{(k-1)k(6n+8k-1)}{6}\\[5pt]
&=|\Pi'|+N\\[5pt]
&=|\pi|+N.
\end{align*}
So $\varphi$ is a bijection satisfying \eqref{fina}, and the proof is complete.
\qed

Combining Theorem \ref{right}, relation \eqref{RRe} and Theorem \ref{La}, we are led to \eqref{mainresult}. This completes  the proof of Conjecture \ref{conjecture}.

\vskip 3mm \noindent {\bf Acknowledgments.}
This work was supported
by the 973 Project, the PCSIRT Project of
the Ministry of Education, the National Science Foundation of
China.

\end{document}